\documentclass[12pt]{amsart}
\usepackage{amssymb}
\usepackage{graphicx}
\input epsf 

\newcommand {\supplus}{\mathop{{\supset}\llap{\raise 
0.5pt\hbox{\normalfont\small+}\hskip 0.5pt}}} 

\newcommand {\subplus}{\mathop{{\subset}\llap{\raise 
0.5pt\hbox{\normalfont\small+}\hskip 0.5pt}}}  

\newcommand {\Cee}    {{\mathbb  C}}

\newcommand {\Nee}    {{\mathbb  N}}
\newcommand {\Mee}    {{\mathbb  M}}

\newcommand {\Zee}    {{\mathbb  Z}}

\newcommand {\fA}     {{\mathfrak{A}}}

\newcommand {\fg}     {{\mathfrak{g}}}    %
\newcommand {\fgl}    {{\mathfrak{gl}}}  %

\newcommand {\fosp}   {{\mathfrak{osp}}}
\newcommand {\fpe}    {{\mathfrak{pe}}}   %

\newcommand {\fq}     {{\mathfrak{q}}}

\newcommand {\fs}     {{\mathfrak{s}}}
\newcommand {\fS}     {{\mathfrak{S}}}

\newcommand {\fsl}    {{\mathfrak{sl}}}

\newcommand {\fspe}   {{\mathfrak{spe}}}

\newcommand {\cal} {\mathcal}

\newcommand {\cL}     {{\cal L}}

\newcommand {\cT}     {{\cal T}}

\def \opname#1#2%
  {\expandafter\newcommand \csname #1\endcsname {{\mathop{#2}\nolimits}}}

\newcommand{\rmname}[1]
  {\expandafter\newcommand \csname #1\endcsname {{\operatorname{#1}}}}

\newcommand{\rmnameii}[2]
  {\expandafter\newcommand \csname #1\endcsname {{\operatorname{#2}}}}

\rmname{act}
\rmname{Ad}
\rmname{Add}
\rmname{ad}
\rmname{Alt}
\rmname{alt}
\rmname{Ann}
\rmname{antidiag}
\rmname{Ber}
\rmname{ber}
\rmname{Br}
\rmname{card}
\rmname{ch}
\rmname{Char}
\rmname{cem}
\rmname{cj}
\rmname{Cliff}
\rmname{cntr}
\rmname{codim}
\rmname{coind}
\rmname{const}
\rmname{col}
\rmname{cork}
\rmname{cpr}
\rmname{diag}
\rmnameii{Div}{div}
\rmname{Def}
\rmname{Der}
\rmname{Dim}
\rmname{End}
\rmname{Even}
\rmname{Ext}
\rmname{gr}
\rmname{Hom}
\rmname{HT}
\rmnameii{Ht}{ht}
\rmname{hwt}
\rmname{Id}
\rmname{id}
\rmname{ind}
\rmname{Ind}
\rmname{Inf}
\rmname{irr}
\rmname{Le}
\rmname{Lie}
\rmname{lwt}
\rmname{mult}
\rmname{Mor}
\rmname{nm}
\rmname{Ob}
\rmname{Odd}
\rmname{Osc}
\rmname{per}
\rmname{Pic}
\rmname{pr}
\rmname{pro}
\rmname{Prime}
\rmname{Proj}
\rmname{prt}
\rmname{pt}
\rmname{Q}
\rmname{qet}
\rmname{qtr}
\rmname{rd}
\rmname{rk}
\rmname{row}
\rmname{Res}
\rmname{salt}
\rmname{Sch}
\rmname{sch}
\rmname{SBr}
\rmname{scalar}
\rmname{Ser}
\rmname{sign}
\rmname{Smbl}
\rmname{spin}
\rmname{ssym}
\rmname{str}
\rmname{st}
\rmname{sgn}
\rmname{sq}
\rmname{symm}
\rmname{supp}
\rmname{Supp}
\rmname{St}
\rmname{Spec}
\rmname{Spm}
\rmname{tr}
\rmname{vpt}
\rmname{weyl}
\rmname{Weyl}
\rmname{Witt}

\opname{vvol}  {{v\hspace{-0.1ex}o\hspace{-0.02ex}l\/}}
\opname{pnt}  {\text{\normalfont pt}}
\opname{Span} {{Span}}
\opname{slim} {\overline{\lim}}
\opname{Vol}  {{V\hspace{-0.55ex}o\hspace{-0.02ex}l\/}}
\opname{Par}  {{P\hspace{-0.3ex}a\hspace{-0.05ex}r\/}}

\rmname{Mat}
\rmname{Bil}
\rmname{Diff}
\rmname{Ker}
\rmname{Herm}
\rmname{Coker}
\rmname{Conn}
\rmname{Covect}
\rmname{Vect}
\rmname{Int}

\rmnameii {IM} {Im}
\rmnameii {RE} {Re}

\opname{Aut} {{A\hspace{-0.2ex}u\hspace{-0.1ex}t\/}}
\opname{GL} {{G\hspace{-0.3ex}L}}
\opname{SL} {{S\hspace{-0.3ex}L}}
\opname{Exp} {{E\hspace{-0.2ex}x\hspace{-0.1ex}p\/}}
\opname{GQ} {{G\hspace{-0.2ex}Q}}
\opname{OSp} {{O\hspace{-0.25ex}S\hspace{-0.15ex}p\/}}
\opname{Out} {{O\hspace{-0.25ex}u\hspace{-0.15ex}t\/}}
\opname{Spp} {{S\hspace{-0.2ex}p\/}}
\opname{SpO} {{S\hspace{-0.2ex}p\hspace{-0.02ex}O\/}}
\opname{Pe} {{P\hspace{-0.25ex}e\/}}
\opname{SPe} {{S\hspace{-0.25ex}P\hspace{-0.25ex}e\/}}
\opname{Spin} {{S\hspace{-0.25ex}p\hspace{-0.05ex}i\hspace{-0.1ex}n\/}}
\opname{Iso} {{I\hspace{-0.25ex}s\hspace{-0.1ex}o\/}}
\opname{SSPe} {{S\hspace{-0.25ex}S\hspace{-0.15ex}P\hspace{-0.25ex}e\/}}
\opname{PeU} {{P\hspace{-0.25ex}e\hspace{-0.1ex}U\/}}
\opname{QU} {{Q\hspace{-0.15ex}U\/}}
\opname{U} {{U\/}}

\opname{cGQ} {{\cal G \hspace{-0.2em} Q \/}}
\opname{cSL} {{\cal S \hspace{-0.2em} L \/}}
\opname{cGL} {{\cal G \hspace{-0.2em} L \/}}
\opname{cOSp} {{\cal O \hspace{-0.2em} S \hspace{-0.3em} \it p\/}}
\opname{cPe} {{\cal P \hspace{-1.5pt} \it e\/}}
\opname{cVect} {{\cal V \hspace{-1.5pt} \it e\hspace{-0.1ex}c\hspace{-0.1ex}t\/}}
\opname{cVol} {{\cal V \hspace{-1.5pt} \it o\hspace{-0.1ex}l\/}}
\opname{cAut} {{\cal A \hspace{-0.2em} \it u\hspace{-0.1em}t\/}}
\opname{cCovect} {{\cal C \hspace{-1.5pt}
     \it o\hspace{-0.1ex}v\hspace{-0.1ex}e\hspace{-0.1ex}c\hspace{-0.1ex}t\/}}
\opname{CW} {{C\hspace{-0.15ex}W}}

\newcommand {\ev} {{\bar0}}
\newcommand {\od} {{\bar1}}
\newcommand {\eps} {\varepsilon}

\newcommand {\tto} {\longrightarrow}

\newcommand {\bcdot}   {\mathbin{\hbox{\raise.4ex\hbox{\bf.}}}} 

\newcommand {\secno} {}
\newcommand {\ssecfont} {\normalfont\bf}

\newtheorem{Theorem}{\secno Theorem}

\newtheorem{Lemma}[Theorem]{\secno Lemma}

\newtheorem{Corollary}[Theorem]{\secno Corollary}

\newenvironment {th*}[1]
    {\gdef\thname{#1} \begin{thn}}%
    {\end{thn}}
\newtheorem{thn}[Theorem] {\thname}

\theoremstyle{definition}

\newenvironment {ex*}[1]
    {\gdef\thname{#1} \begin{exn}}%
    {\end{exn}}
\newtheorem{exn}[Theorem]{\thname}

\theoremstyle{remark}

\newtheorem{Remark}[Theorem]{\secno Remark}

\newenvironment {rem*}[1]
    {\gdef\thname{#1} \begin{remn}}%
    {\end{remn}}
\newtheorem{remn}[Theorem]{\thname}

\newcommand {\ssec}{\subsection*}

\newcommand {\ssbegin}[2]
  {\def \secno {\gdef \secno {}{\ssecfont #1. }}%
   \begin{#2}}
\begin{document}
\title[Classical invariant theory for Lie superalgebras]{An analog of 
the classical invariant theory \\ for Lie superalgebras. II}

\author{Alexander Sergeev}

\address{Dept.  of Math., Univ.  of Stockholm, Roslagsv.  101, 
Kr\"aftriket hus 6, S-106 91, Stockholm, Sweden (On leave of absence 
from Balakovo Inst.  of Technology Technique and Control)\\ e-mail: 
mleites@matematik.su.se subject: for Sergeev}

\thanks{I am thankful to D. Leites for support and help.}

\begin{abstract} Let $V$ be a finite-dimensional superspace over 
$\Cee$ and $\fg$ a simple (or a ``close'' to simple) matrix Lie 
superalgebra, i.e., a Lie subsuperalgebra in $\fgl(V)$.  Under the 
{\it classical invariant theory} for $\fg$ we mean the description of 
$\fg$-invariant elements of the algebra
$$
\fA^{p, q}_{k, l}=S^{\bcdot}(V^{k}\oplus \Pi (V)^{l}\oplus V^{*p}\oplus 
\Pi (V)^{*q}),
$$
where $V^p$ denotes $V\oplus \dots\oplus V$ ($p$ summands).  We give 
such description for $\fgl(V)$, $\fsl(V)$ as well as $\fosp(V)$ and 
their ``odd'' analogs: $\fq(V)$, $\fs\fq(V)$; $\fpe(V)$ and 
$\fspe(V)$.

In \cite{S2} the description of $\fg$-invariant elements in $\fA^{p, 
q}_{k, l}$ was given up to polarization operators, i.e., the elements 
of $U(\fgl (U\otimes W))$ which naturally act on $\fA^{p, q}_{k, l}$ 
presented in the form $S^{\bcdot}(U\otimes V\bigoplus V^{*}\otimes 
W)$.  Here we give a complete description of the generators in the 
algebra of invariants and describe the relations between the 
invariants of the scalar product type.

\end{abstract}

\subjclass{17A70 (Primary) 13A50 (Secondary)}

\keywords{Invariant theory, Lie superalgebras.}

\maketitle

This paper is a detailed exposition of \cite{S3} with several new 
results added. It also complements and refines the results of \cite{S2}.

\section*{\S 1. Preliminaries}

In what follows $\fS_{k}$ stands for the symmetric group on $k$ 
elements. Let $\lambda$ be a partition of the number $k$ and $t$ a 
$\lambda$-tableau. Recall that $t$ is called {\it standard} if the 
numbers in its rows and columns grow from left to right and 
downwards. Denote by $C_{t}$ the column stabilizer of $t$, ler $R_{t}$ 
be its row stabilizer. We further set
$$
e_{t}=\mathop{\sum}\limits_{\tau\in C_{t};\;\sigma\in 
R_{t}}\eps(\tau)\sigma\tau, \quad 
\tilde e_{t}=\mathop{\sum}\limits_{\tau\in C_{t};\;\sigma\in 
R_{t}}\eps(\tau)\tau\sigma.\eqno{(0.1)}
$$

Let $\Nee$ be the set of positive integers, $\bar\Nee$, another ,
``odd'', copy of $\Nee$ and let $\Mee=\Nee\coprod\bar\Nee$ be ordered so 
that each element of the ``even'' copy, $\Nee$, is smaller than any 
element form the ``odd'' copy, while inside of each copy the order is 
the natural one. We will call the elemnts from $\Nee$ ``even'' and 
those form $\bar\Nee$ ``odd'' ones; so we can encounter an ``even'' 
odd element, etc.

Let $I$ be the sequence of elements from $\Mee$ of length $k$.  Let us 
fill in the tableau $t$ with elements from $I$ replacing element 
$\alpha$ with $i_{\alpha}$.  The sequence $I$ is called $t$-{\it 
semistandard} if the elements of $t$ do not decrease from left to 
right and downwards; the ``even'' elements strictly increase along 
columns; the ``odd'' elements strictly increase along rows.

The group $\fS_{k}$   naturally acts on sequences $I$. Let $\fA$ be 
the free supercommutative superalgebra with unit generated by 
$\{x_{i}\}_{i\in I}$. For any $\sigma\in\fS_{k}$ define $c(I, 
\sigma)=\pm 1$ from the equation
$$
c(I, \sigma)x_{I}=x_{\sigma^{-1}I}\text{ where }x_{I}=x_{i_{1}}\dots 
x_{i_{k}}.\eqno{(0.2)}
$$
Clearly, $c(I, \sigma)$ is a cocycle, i.e.,
$$
c(I, \sigma\tau)=c(\sigma^{-1}I, \tau)c(I, \sigma).
$$
With the help of this cocycle a representation of $\fS_{k}$ in 
$T^{k}(V)=V^{\otimes k}$ for any
superspace $V$ is defined:
$$
\sigma v_{I}= c(I, \sigma^{-1})v_{\sigma I}, \text{ where $v_{I}= 
v_{i_{1}}\otimes \dots \otimes v_{i_{k}}$ and $v_{i_{\alpha}}\in V$ for each 
$\alpha$}.\eqno{(0.3)}
$$
Let $\{v_{1}, \dots , v_{n}; v_{\bar 1}, \dots , v_{\bar  m}\}$ be a basis 
of $V$ in the standard format (the even elements come first followed 
by the odd ones). Then the elements $v_{I}$ for all possible 
sequences $I$ of length $k$ and with elements from 
$$
R_V=\{1, \dots , n; \bar 
1, \dots , \bar  m\}\eqno{(0.4)}
$$ 
form a basis of $T^{k}(V)$.

The following theorem describes the decomposition of $T^{k}(V)$ into 
irreducible $\fS_{k}\times \fgl(V)$-modules.

\ssbegin{1.1}{Theorem} {\em (Cf. \cite{S1}.)} The commutant of the 
natural $\fgl(V)$-action on $T^{k}(V)$ is isomorphic to $\Cee[\fS_{k}]$ 
and
$$
T^{k}(V)=\mathop{\oplus}\limits_{\lambda: \lambda_{n+1}\leq 
m}S^{\lambda}\otimes V^{\lambda},
$$
where $S^{\lambda}$ is an irreducible $\fS_{k}$-module and $V^{\lambda}$ 
is an irreducible $\fgl(V)$-module.
\end{Theorem}

The following refinement of Theorem 1.1 holds:

\ssbegin{1.2}{Theorem} If $t$ runs over the standard tableaux of type 
$\lambda$ and $I$ runs over semistandard $t$-sequences, then the 
family $\{e_{t}(v_{I})\}$ (resp.  $\{\tilde e_{t}(v_{I})\}$) is a 
basis in $S^{\lambda}\otimes V^{\lambda}$.  Moreover, for a fixed $t$ 
the families $\{e_{t}(v_{I})\}$ and $\{\tilde e_{t}(v_{I})\}$ span 
$V^{\lambda}$.
\end{Theorem}

Proof follows from results of \cite{S1}. \qed

Let $U$ and $W$ be two superspaces with bases $u_{i}$ and $w_{j}$, 
where $i\in R_U$, $j\in R_W$,
respectively, and where 
$$
R_U=\{1, \dots , k; \bar 1, \dots , \bar l\}\text{ and }
R_W=\{1, \dots , p; \bar 1, \dots , \bar q\}.
$$
The symmetric algebra $S^{\bcdot}(U\otimes W)$ is generated by 
$z_{ij}=u_{i}\otimes w_{j}$ for $i\in R_U$ and $j\in R_W$.  Let $I$ be a 
sequence of length $N$ with elements from $R_U$ and $J$ a sequence of 
the same length with elements form $R_W$.  Let $p(i_{\alpha})$ and 
$p(j_{\beta})$be the parities of the corresponding elements of the 
sequence.  Set $\alpha(I, J)=\mathop{\sum}\limits_{\alpha 
>\beta}p(i_{\alpha})p(j_{\beta})$ and define an element of $S^{\bcdot} (U\otimes 
W)$ by setting
$$
Z(I, J)=(-1)^{\alpha(I, 
J)}\mathop{\prod}\limits_{\alpha=1}^{N}Z_{i_{\alpha}j_{\beta}}.
\eqno{(1.1)}
$$
For a given tableau $t$ of order $N$ we define polynomials
$$
\renewcommand{\arraystretch}{1.4}
\begin{array}{l}
    P_{t}(I, J)=\mathop{\sum}\limits_{\sigma\in R_{t}, \, \tau\in 
    C_{t}}\eps(\tau)c(I, (\sigma\tau)^{-1})Z(\sigma\tau I, J),\\
\tilde P_{t}(I, J)=\mathop{\sum}\limits_{\sigma\in R_{t}, \, \tau\in 
C_{t}}\eps(\tau)c(I, (\tau\sigma)^{-1})Z(\tau\sigma I, J).
\end{array}
$$
The Lie superalgebras $\fgl(U)$ and $\fgl(W)$ naturally act on 
$S^{\bcdot}(U\otimes W)$ and their actions commute.

\ssbegin{1.3}{Theorem} $S^{\bcdot}(U\otimes 
W)=\mathop{\oplus}\limits_{\lambda}U^{\lambda}\otimes W^{\lambda}$, 
where $U^{\lambda}$ and $W^{\lambda}$ are irreducible $\fgl(U)$- and 
$\fgl(W)$-modules, respectively, corresponding to the partition 
$\lambda$ and the sum runs over partitions such that $\lambda_{\alpha 
+1}\leq \beta$ for $\alpha=\min(k, p)$ and $\beta=\min(l, q)$.
\end{Theorem}

\begin{proof} By Theorem 1.1
$$
W^{\otimes N}=\oplus ~  W^{\lambda}\otimes S^{\lambda}\text{ and }
U^{\otimes N}=\oplus ~U^{\mu}\otimes S^{\mu}.
$$
Hence,
$$
\renewcommand{\arraystretch}{1.4}
\begin{array}{l}
    S^{N}(U\otimes W)=((U\otimes W)^{\otimes N})^{\fS_{N}}=
(U^{\otimes N}\otimes W^{\otimes N})^{\fS_{N}}=\\
\mathop{\oplus}\limits_{\lambda, \mu} (U^{\lambda}\otimes 
W^{\mu}\otimes S^{\lambda}\otimes S^{\mu})^{\fS_{N}}= 
\mathop{\oplus}\limits_{\lambda, \mu} (U^{\lambda}\otimes 
W^{\mu})\otimes (S^{\lambda}\otimes S^{\mu})^{\fS_{N}}.
\end{array}
$$
Since $(S^{\lambda})^{*}\simeq S^{\lambda}$ and since $S^{\lambda}$ 
and $S^{\mu}$ are irreducible, we have
$$
(S^{\lambda}\otimes S^{\mu})^{\fS_{N}}=\Hom_{\fS_{N}}(S^{\lambda}, 
S^{\mu})=0 \text{ if $\lambda\neq \mu$ and $\Cee$ otherwise}.
$$
Theorem is proved. \end{proof}

\ssbegin{1.4}{Theorem} Let $t$ be a standard tableau of type $\lambda$ 
and let $I$ and $J$ be $t$-semistandard sequences.  Then the family 
$P_{t}(I, J)$, as well as the similar family 
$\tilde P_{t}(I, J)$, forms a basis in the module $U^{\lambda}\otimes 
W^{\lambda}$.
\end{Theorem}

\begin{proof} The natural homomorphism
$$
\phi_{N}: U^{\otimes N}\otimes W^{\otimes N}\tto S^{N}(U\otimes W)
$$
is, clearly, a homomorphism of $\fgl(U)\oplus \fgl(W)$-modules. It is 
not difficult to verify that
$$
\phi_{N}(e_{t}(v_{I})\otimes\tilde e_{t}(w_{J}))=c\cdot P_{t}(I, 
J)\text{ for a constant } c.
$$

Let $t$ be a fixed $\lambda$-tableau, $I$ and $J$ two $t$-semistandard 
sequences with elements from $R_U$ and $R_W$, respectively.  Then by 
Theorem 1.2 the vectors $e_{t}(v_{I})\otimes\tilde e_{t}(w_{J})$ form 
a basis of a subspace $L\subset U^{\otimes N}\otimes W^{\otimes N}$ 
which is also a $\fgl(U)\oplus \fgl(W)$-submodule.  By the same 
theorem, $L\simeq U^{\lambda}\otimes W^{\lambda}$ and it remains to 
show that $\phi_{N}(L)\neq 0$. For this it suffices to show that 
there exists an $l\in L$ such that $\phi(l)\neq 0$. Since 
$\phi_{N}(\sigma v_{i}\otimes \sigma w_{J})=
\phi_{N}(v_{I}\otimes  w_{J})$, it follows that
$$
\renewcommand{\arraystretch}{1.4}
\begin{array}{l}
    \phi_{N}(e_{t}(v_{I})\otimes \tilde e_{t}(w_{J}))=
c\phi_{N}(e_{t}(v_{I})\otimes w_{J})=\\
c\phi_{N}(\sigma e_{t}(v_{I})\otimes \sigma w_{J})=
c\phi_{N}(e_{\sigma t}(\sigma v_{I})\otimes \sigma w_{J})=\\
\pm c\phi_{N}(e_{\sigma t}(v_{\sigma I})\otimes w_{\sigma J}).
\end{array}
$$
Therefore, we may assume that the tableau $t$ is consequtively filled 
in along the rows with the numbers 1, 2, etc.  Observe that the 
sequences $\sigma I$ and $\sigma J$ remain $\sigma t$-semistandard.

Let $I=J$ be the sequence
$$
\underbrace{1\dots 1}_{\lambda_{1}}
\underbrace{2\dots 2}_{\lambda_{2}}\dots
\underbrace{\alpha\dots \alpha}_{\lambda_{\alpha}}\bar 1\dots 
\bar\lambda_{\alpha+1}\bar 1\dots 
\bar\lambda_{\alpha+2}\dots \bar 1\dots 
\bar\lambda_{\gamma},
$$
where $(\lambda_{1}, \lambda_{2}, \dots, \lambda_{\alpha}, \dots , 
\lambda_{\gamma})$ is the partition corresponding to $t$, 
$\alpha=\min(\dim~U_{\ev}, \dim~W_{\ev})$, 
$\beta=\min(\dim~U_{\od}, \dim~W_{\od})$.

It is not difficult to verify that $\phi_{N}(e_{t}(v_{I})\otimes 
w_{I})\neq 0$.

Since $\phi_{N}$ is a homomorphism of $\fgl(U)\oplus \fgl(W)$-modules, 
its restriction onto $L$ is an isomorphism.  This implies the 
statement of Theorem for the family $P_{t}(I, J)$.  For the family 
$\tilde P_{t}(I, J)$ proof is similar.
\end{proof}

Let us elucidate how the results obtained can be applied to the 
invariant theory.

Let $\fg\subset \fgl(V)$ be a Lie superalgebra. Under ``the invariant 
theory of $\fg$'' we understand the description of $\fg$-invariants in 
the superalgebra $\fA^{p, q}_{k, l}=S(U\otimes V\; \oplus \;
V^{*}\otimes W)$.

On $\fA^{p, q}_{k, l}$, the Lie superalgebras $\fgl(U)$ and $\fgl(W)$ 
naturally act. By Theorem 1.3 we have
$$
\fA^{p, q}_{k, l}=\mathop{\oplus}\limits_{\lambda, \mu} 
\; U^{\lambda}\otimes V^{\lambda}\otimes V^*{}^{\mu}\otimes W^{\mu}.
$$
Therefore, to describe $\fg$-invariant elements, it suffices to 
describe the $\fg$-invariants in $V^{\lambda}\otimes 
V^*{}^{\mu}=\Hom(V^{\mu}, V^{\lambda})$.  But $(V^{\lambda}\otimes 
V^*{}^{\mu})^\fg=\Hom_{\fg}(V^{\mu}, V^{\lambda})$, i.e., the 
description of $\fg$-invariants is equivalent to the description of 
$\fg$-homomorphisms of $\fg$-modules $V^{\mu}$.

Let us consider how the method works in the simplest example: 
$\fg=\fgl(V)$.  Let $\{e_{i}: i\in R_V\}$ be a basis of $V$ in a 
standard format; $\{e_{i}^{*}\}$ the left dual basis.  Set
$$
\theta=\mathop{\sum}\limits_{i\in T}e_{i}\otimes e_{i}^{*}, \quad 
\hat\theta=\mathop{\sum}\limits_{i\in T}(-1)^{p(i)}e_{i}^{*}\otimes e_{i}. 
$$
It is not difficult to verify that $\theta$ and $\hat\theta$ are
$\fg$-invariants.

Set 
$$
T^{p, q}(V)=V^{\otimes p}\otimes V^*{}^{\otimes q}, \quad 
\hat T^{p, q}(V)=V^*{}^{\otimes p}\otimes V^{\otimes q}.  
$$
On $T^{p, q}(V)$ and $\hat T^{p, q}(V)$, the group $\fS_{p}\times 
\fS_{q}$ acts and its action commutes with that fo $\fgl(V)$.  Hence, 
$\fS_{p}\times \fS_{q}$ also acts on the space of $\fgl(V)$-invariants 
in $T^{p, q}(V)$ and $\hat T^{p, q}(V)$.

\section*{\S 2. The invariants of \protect $\fgl(V)$}

Set $v_r{}^*=(x_{r1}, \dots , x_{rn}; x_{r\bar 1}, \dots , x_{r\bar 
m})$ and $v_s=(x_{1s}{}^*, \dots , x_{n s}{}^*; x_{\bar 1s}{}^*, \dots 
, x_{\bar m s}{}^*)^{t}$, where $x_{ri}=u_r\otimes e_i$ and 
$x_{is}=e_i{}^*\otimes w_s$, i.e., $v_r{}^*$ is a row vector and $v_s$ 
is a column vector, so their scalar product is equal to $(v_r{}^*, 
v_s)=\mathop{\sum}\limits_i x_{ri}x_{i s}{}^*$.

\ssbegin{2.1}{Theorem} The algebra of $\fgl(V)$-invariant
elements in $\fA^{p, q}_{k ,l}$ is generated by the elements 
$(v_r{}^*, v_s)$ for all $r\in R_U$, $s\in R_W$.
\end{Theorem}

\begin{proof} Let $A$ be a supercommutative superalgebra,
$L$ a $\fg$-module; let  $L_A=(L\otimes A)_{\ev}$ and
$\fg_A=(\fg\otimes A)_{\ev}$. 

The elements of $S^{\bcdot}(L^*)$ may be considered as functions on
$L_A$ with values in $A$. Let $l\in L_A=(L\otimes
A)_{\ev}=(\Hom(L^*,  A))_{\ev}$. Therefore, $l$ determines a
homomorphism $\phi_l: S(L^*)\tto A$. Set
$$
f(l)=\phi_l(f)\text{ for any }f\in S^{\bcdot}(L^*).
$$
Observe that $\fg_A$ naturally acts on $L_A$ and on the
algebra of functions on $L_A$. \end{proof}

\ssbegin{2.1.1}{Statement} {\em (\cite{S2})} Let $A$ be the
Grassmann superalgebra with more indeterminates than $\dim
L_\od$. The element $x$ of $S(L^*)$ is a $\fg$-invariant if and
only if $x$ is a $\fg_A$-invariant when considered as a
function on $L_A$.
\end{Statement}

\ssec{2.1.2.  Proof of Theorem 2.1} Let $V^p$ denote $V\oplus 
\dots\oplus V$ ($p$ summands).  Set $L=V^p\oplus \Pi(V)^q\oplus 
(V^*)^k\oplus \Pi(V^*)^l$; then $S(L^*)=\fA^{p, q}_{k, l}$ and we can 
consider $L_A$ as the set of collections
$$
\cL=(v_1, \dots , v_p, v_{\bar 1}, \dots , v_{\bar q}, v_1^*,
\dots , v_k^*, v_{\bar 1}^*, \dots , v_{\bar l}^*),
$$
where $v_s\in V\otimes A$ and $v_t^*\in
\Hom_A(V\otimes A,  A)$ and the parity of these
vectors coincide with the parities of their indices.

Let us write the vectors with right coordinates and the
covectors with left ones:
$$
v_s=\mathop{\sum}\limits_i e_ia_{is}^*, \quad 
v_t^*=\mathop{\sum}\limits_i a_{ti}e_i^*.
$$
Consider now the elements of $\fA^{p, q}_{k, l}$ as functions on
$\cL$, by setting 
$$
x_{is}^*(\cL)=a_{is}^*, \quad x_{ti}^*(\cL)=a_{ti}.
$$
Therefore, thanks to 2.1.1 it suffices to describe the functions on $\cL$
contained in the subalgebra generated by the coordinate
functions {\bf ??} and invariants with respect to $GL(V\otimes A)$.

Now, since the scalar products turn into scalar products under
the $\fgl(U)\oplus  \fgl(W)$-action, it is sufficient to confine
ourselves to the invariants in $\fA^{n, m}_{n, m}$.

Denote by $M$ the set of collections $(v_1, \dots , v_n,
v_{\bar 1}, \dots , v_{\bar m})$ that form bases of $V\otimes
A$. In Zariski topology the set $M$ is dense in the space of all
collections. If $f$ is an invariant and $\cL\in M$, then there
exists $g\in GL(V\otimes A)$ such that $gv_i=e_i$ for each $i\in
T$.

Therefore, $f(\cL)=f(g\cL)=f(e_1, \dots , e_n,
gv_{\bar 1}^*, \dots , gv_{\bar m}^*)$ and $f(\cL)$ is a
polynomial in coordinates of the $gv_{\bar t}$. But $(gv_{\bar
t}^*, e_i)=(v_{\bar t}^*, g^{-1}e_i)=(v_{\bar t}^*, v_i)$ which
proves the theorem. \qed

\begin{Corollary} The nonzero $\fgl(V)$-invariants in $T^{p, q}$
only exist if $p=q$. In this case the $\fS_p\times
\fS_p$-module of invariants is generated by the 
images of the canonical elements $\theta^{\otimes p}$ in $T^{p, p}$ and 
$\hat\theta$ in $\hat T^{p, p}$.
\end{Corollary}

Consider now the algebra homomorphism
$$
S^{\bcdot}(U\otimes W)\tto (\fA^{p, q}_{k, l})^{\fgl(V)}, \quad 
u_{r}\otimes w_s\mapsto (v_r^*, v_s).\eqno{(2.1)}
$$
The kernel of this homomorphism is the ideal of relations
between the scalar products.

\ssbegin{2.2}{Theorem} The ideal of relations between scalar products 
$(v_r^*, v_s)$ is generated by the polynomials $P_t(I, J)$, where $t$ is a 
fixed standard rectangular $(n+1)\times (m+1)$ tableau, $I$ and $J$ 
are $t$-semistandard sequences with elements from $R_U$ and $R_W$, 
respectively.  \end{Theorem}

\begin{proof} By Theorem 1.3 $S^{\bcdot}(U\otimes 
W)=\mathop{\oplus}\limits_{\lambda}U^\lambda\otimes W^\lambda$ and
$$
\fA^{p, q}_{k, l}=S^{\bcdot}(U\otimes V\oplus V^*\otimes 
W)=S^{\bcdot}(U\otimes V)\otimes S^{\bcdot}(V^*\otimes 
W)=(\mathop{\oplus}\limits_{\mu}U^\mu\otimes V^\mu)\otimes 
(\mathop{\oplus}\limits_{\nu}(V^*)^\nu\otimes W^\nu),
$$
hence, $(\fA^{p, q}_{k, l})^{\fgl(V)}=\mathop{\oplus}\limits_{\mu: 
\mu_{n+1}\leq m}U^\mu\otimes W^\mu$.  Since homomorphism (2.1)  is a 
homomorphism of $\fgl(V)\oplus  \fgl(W)$-modules, its kernel coinsides 
with $\mathop{\oplus}\limits_{\lambda: \lambda_{n+1}\geq 
m+1}U^\lambda\otimes W^\lambda$.

Let $\nu$ be a $(n+1)\times (m+1)$ rectangle.  The condition 
$\lambda_{n+1}\geq m+1$ means that $\lambda\supset \nu$ and by Theorem 
1.3  it suffices to demonstrate that $P_t(I, J)$, where $t$ is a fixed 
standard rectangular tableau of size $\lambda$, belongs to the ideal 
generated by $U^\nu\otimes W^\nu$.

Let $e_t$ be the corresponding minimal idempotent, $e_s$ the minimal 
idempotent for a standard tableau $s$ of size $\nu$.  Decomposing 
$R_t$ with respect to the right cosets relative $R_s$ and decomposing 
$C_t$ with respect to the left cosets relative $C_s$ we obtain a 
representation of $e_t$ in the form $\sum\tau_ie_s\sigma_j$.  This 
implies that $P_t(I, J)$ is the sum of polynomials of the form 
$f_iP_{t_{i}}(I_{i}, J_{j})\phi_{j}$, i.e., belongs to the ideal 
generated by the $P_t(I, J)$.
\end{proof}

\section*{\S 3. The invariants of \protect $\fsl(V)$}

First, let us describe certain tensor invariants.  Obviously, all 
$\fgl(V)$-invariants are also $\fsl(V)$-invariants.  Therefore, we 
will only describe the $\fsl(V)$-invariants which are not 
$\fgl(V)$-invariants.  Denote by $\theta_k=\theta^{\otimes k}$ the 
invariant in $T^{k, k}$ and by $\hat\theta_k=\hat\theta^{\otimes k}$ the 
invariant in $\hat T^{k, k}$, and for a given sequence $I$ with elements 
from $R_V$ set
$$
v_I=e_{i_{1}}\otimes \dots \otimes e_{i_{k}}\text{ and }
v_I^*=e_{i_{1}}^*\otimes \dots \otimes e_{i_{k}^*}.
$$
Let us represent $\fsl(V)$ in the form 
$\fg=\fg_-\oplus\fg_0\oplus\fg_+$, where $\fg_0=\fg_\ev$ and $\fg_\pm$ 
are the $\fg_0$-modules generated by the positive and negative root 
vectors, respectively.

Let $\{X_\alpha\}_{\alpha\in R^-}$ and $\{X_\beta\}_{\beta\in R^+}$, 
where $R^{\pm}$ are the sets of positive (negative) roots, be some bases of $\fg_-$ and 
$\fg_+$, respectively; set $X_-=\prod X_\alpha$ and $X_+=\prod 
X_\beta$.  The elements $X_\pm$ are uniquely determined up to a 
constant factor because the subalgebras $\fg_\pm$ are commutative.

\ssbegin{3.1}{Lemma} Let $M$ be a $\fg_0$-module and $\tilde
M=\ind^{\fg}_{\fg_0}(M)$ be the induced $\fg$-module. Then
each of the correspondences $m\mapsto X_+X_-m$ and $m\mapsto
X_-X_+m$ is a bijection of $M^{\fg_{0}}$ onto $\tilde M^{\fg}$.
\end{Lemma}

\begin{proof} As follows from Lemma 4.2 below, $\dim M^{\fg_{0}}=\dim 
\tilde M^{\fg}$.  Therefore, it suffices to show that the 
correspondence $m\mapsto n=X_+X_-m$ is injective map of $M^{\fg_{0}}$ 
to $\tilde M^{\fg}$.  The injectivity is manifest, so we only have to 
check that the image is $\fg$-invariant.  Clearly, $\fg_+n=\fg_0n=0$.  
Therefore, it suffices to verify that $X_{-\alpha}n=0$ for every 
simple root $\alpha$.  This is subject to a direct check with the help 
of the multiplication table in $\fgl(V)$.  \end{proof}

\ssbegin{3.2}{Lemma} Let $V_1$ and $V_2$ be finite
dimensional  $\fg_0$-modules. Set  $\fg_+V_1=\fg_-V_2=0$.
Then 
$$
\ind^{\fg}_{\fg_0\oplus\fg_+}V_1
\otimes\ind^{\fg}_{\fg_0\oplus\fg_-}(V_2)
\simeq \ind^{\fg}_{\fg_0}(V_1\otimes V_2).\eqno{(3.2.1)}
$$
is an isomorphism of $\fg$-modules. \end{Lemma}

\begin{proof} Since the dimensions of both modules are equal,
it suffices to show that the natural homomorphism
$$
\ind^{\fg}_{\fg_0}(V_1\otimes
V_2)\tto \ind^{\fg}_{\fg_0\oplus\fg_+}V_1
\otimes\ind^{\fg}_{\fg_0\oplus\fg_-}(V_2).\eqno{(3.2.2)} 
$$
is surjective, i.e., the module generated by $V_1\otimes V_2$
coinsides with the whole module.  

The module in the right hand side has a natural filtration induced by 
filtrations of the modules $\ind^{\fg}_{\fg_0\oplus\fg_+}V_1$ and 
$\ind^{\fg}_{\fg_0\oplus\fg_-}(V_2)$.  Let the $X_\alpha$ be a basisi 
of $\fg_+$ and $X_{-\alpha}$ a basisi of $\fg_-$.  Consider the module 
$W$ generated by $V_1\otimes V_2$, i.e., by the elements of filtration 
zero and the element
$$
w=X_{-\alpha_{1}}\dots X_{-\alpha_{k}}v_1\otimes 
X_{\beta_{1}}\dots X_{\beta_{l}}v_2.
$$
We have
$$
\renewcommand{\arraystretch}{1.4}
\begin{array}{l}
    w=X_{-\alpha_{1}}(X_{-\alpha_{2}}\dots
X_{-\alpha_{k}}v_1\otimes  X_{\beta_{1}}\dots
X_{\beta_{l}}v_2)
\pm X_{-\alpha_{2}}\dots
X_{-\alpha_{k}}v_1\otimes  X_{-\alpha_{1}}X_{\beta_{1}}\dots
X_{\beta_{l}}v_2=\\
X_{-\alpha_{1}}(X_{-\alpha_{2}}\dots
X_{-\alpha_{k}}v_1\otimes  X_{\beta_{1}}\dots
X_{\beta_{l}}v_2)\pm \\
\pm X_{-\alpha_{2}}\dots
X_{-\alpha_{k}}v_1\otimes (\mathop{\sum}\limits_i
X_{\beta_{1}}\dots X_{\beta_{i-1}} [X_{-\alpha_{1}},
X_{\beta_{i}}]X_{\beta_{i+1}}\dots X_{\beta_{l}}v_2. 
\end{array}
$$
Since each summand is of filtration $<k+l$, they
belong to $W$ by inductive hypothesis; hence, so
does $w\in W$. \end{proof}

Let $t$ be a tableau consisting of $m$ columns and $n+k$ rows and 
filled in as follows: first, we fill in the tableau $t_1$ that 
occupies the first $n$ rows, next, the tableau $t_2$ that occupies the 
remaining rows, both tableaux are filled in consequtevely column-wise.

Let $s$ be a tableau consisting of $n$ rows and $k+m$ columns and 
filled in as follows: first, we fill in the tableau $s_1$ that 
occupies the first $k$ columns, next, the tableau $t_2$ that occupies 
the remaining columns, both tableaux are filled in consequtevely 
column-wise.

Let $I_k$ be the sequence obtained by $k$-fold repetition of the 
sequence $1, 2, \dots , n$; let $J_k$ be the sequence consisting of 
$k$ copies of $\bar 1$ in a row, next $k$ copies of $\bar 2$ in a row, 
etc., $k$ copies of $\bar m$ in a row.

\ssbegin{3.3}{Theorem} In $\hat T^{(m+k)n, (n+k)m}(V)$, the element
$$
e_s\times \tilde
e_t(v^*_{I_{k}}\otimes\hat\theta_{nm}\otimes v_{J_{k}}) 
$$
is an $\fsl(V)$-invariant.
\end{Theorem}

\begin{proof} Let $N$ be any positive integer.  Consider the map 
$\phi: T^{N, N}(V)\tto T^{N, N}(V_\ev)$ such that 
$\phi(V_\od)=\phi(V^*_\od)=0$.  Clearly, $\phi$ is an $\fS_N\times 
\fS_N$-module homomorphism because it is induced by projections of $V$ 
and $V^*$ onto their even parts.

Take $X_+$ and $X_-$ from Lemma 3.1 and consider the map 
$$
\psi: T^{N, N}(V_\ev)\tto  T^{N,
N}(V),\quad v_0\mapsto X_+X_-v_0.
$$
Clearly, $\psi$ is an $\fS_N\times \fS_N$-module
homomorphism.

Let us consider the restrictions of the maps $\phi$ and $\psi$
onto $T^{N, N}(V)^{\fgl(V)}$ and $T^{N, N}(V_\ev)^{\fgl(V_\ev)}$, 
respectively.

Clearly, $\phi$ sends the first of these spaces into the second one, 
whereas by Lemma 3.1  $\psi$ sends the second of these spaces into the 
first one.  Theorem 1.1  implies that, as $\fS_N\times \fS_N$-modules, 
the spaces $T^{N, N}(V)^{\fgl(V)}$ and $T^{N, N}(V_\ev)^{\fgl(V_\ev)}$ 
have simple spectra.

Let $S^\lambda\otimes S^\lambda\subset T^{N,
N}(V)^{\fgl(V)}$ while $S_0^\lambda\otimes S_0^\lambda\subset
T^{N, N}(V_\ev)^{\fgl(V_\ev)}$ correspond to a typical diagram
$\lambda$ and both are nonzero, i.e., $\lambda_n\geq m$ and 
$\lambda_{n+1} =0$. Then the simplicity of the spectrum and
Lemma 3.1 imply that $\phi$ and $\psi$ are, up to a constant
factor, mutually inverse isomorphisms of the modules
$S^\lambda\otimes S^\lambda$ and $S_0^\lambda\otimes
S_0^\lambda$. 

Let
$$
C_t=\mathop{\cup}\limits_\tau\tau(C_{t_{1}}\times C_{t_{2}}),\quad
R_s=\mathop{\cup}\limits_\sigma\sigma(R_{s_{1}}\times R_{s_{2}})
$$
be the decomposition of the column stabilisor $C_t$ of the tableau
$t$ with respect to the left cosets relative the product of the
column stabilisors of $t_1$ and $t_2$ and  same of the row stabilisor
$R_s$. Then
$$
\tilde e_t=\mathop{\sum}\limits_\tau\eps(\tau)\tau\tilde e_{t_{1}}\tilde
e_{t_{2}},\quad
e_s=\mathop{\sum}\limits_\sigma \sigma e_{s_{1}}e_{s_{2}}.
$$
It is easy to verify that 
$$
\fg_+(V_\ev)=\fg_-(V_\ev ^*)=\fg_+(V_\od
^*)=\fg_-(V_\od)=0.\eqno{(3.3)}
$$

The vector $X_+e_s(v^*_{I_{k}}\otimes v^*_{I_{m}})$ belongs to a 
typical module, is a highest one with respect to 
$\fg_+\oplus(\fg_0)_+$ and nonzero, where $(\fg_0)_+$ is the set of 
strictly upper-triangular matrices with respect to the fixed basis of 
$V_\ev$.  But (0.2) implies that $e_s(v^*_{I_{k}}\otimes v^*_{J_{n}})$ 
is also highest with respect to $\fg_+\oplus(\fg_0)_+$ and lies in the 
same module.  This shows that
$$
X_+e_s(v^*_{I_{k}}\otimes v^*_{I_{m}})=c\cdot 
e_s(v^*_{I_{k}}\otimes v^*_{J_{n}}),\text{ where }c\neq 0.
$$
Further, from Lemmas 3.1 and 3.2 it follows that the vector
$$
X_-X_+\left[e_s(v^*_{I_{k}}\otimes v^*_{I_{m}})\otimes\tilde
e_t(v_{I_{m}}\otimes v_{J_{k}})\right]
$$
is $\fg$-invariant because $e_s(v^*_{I_{k}}\otimes 
v^*_{I_{m}})\otimes\tilde e_t(v_{I_{m}}\otimes v_{J_{k}})$ is 
$\fg_{\ev}$-invariant.  We make use of the fact that $X_+\tilde 
e_t(v^*_{I_{m}}\otimes v^*_{J_{k}})=0$ to decuce that
$$
\renewcommand{\arraystretch}{1.4}
\begin{array}{l}
w=X_-X_+\left[e_s(v^*_{I_{k}}\otimes v^*_{I_{m}})\otimes\tilde
e_t(v_{I_{m}}\otimes v_{J_{k}})\right]=
X_-\left[[X_+e_s(v^*_{I_{k}}\otimes v^*_{I_{m}})]\otimes\tilde
e_t(v_{I_{m}}\otimes v_{J_{k}})\right]=\\
const\cdot X_-[e_s(v^*_{I_{k}}\otimes
v^*_{J_{n}})\otimes\tilde e_t(v_{I_{m}}\otimes v_{J_{k}})]=\\
\mathop{\sum}\limits_{\sigma, \tau}\eps(\tau)\sigma\times
\tau\left (X_-[e_{s_{1}}e_{s_{2}}(v^*_{I_{k}}\otimes
v^*_{I_{m}})\otimes\tilde e_{t_{1}}\tilde
e_{t_{2}}(v_{I_{m}}\otimes v_{J_{k}})]\right)=\\
\mathop{\sum}\limits_{\sigma, \tau}\eps(\tau)\sigma\times
\tau\left (e_{s_{1}}(v^*_{I_{k}})\otimes
X_-[e_{s_{2}}v^*_{I_{m}})\otimes\tilde
e_{t_{1}}(v_{I_{m}})]\tilde e_{t_{2}}(v_{J_{k}})\right)=\\
\mathop{\sum}\limits_{\sigma, \tau}\eps(\tau)\sigma\times
\tau\left (e_{s_{1}}(v^*_{I_{k}})\otimes
X_-X_+[e_{s_{2}}(v^*_{I_{m}})\otimes\tilde
e_{t_{1}}(v_{I_{m}})]\tilde e_{t_{2}}(v_{J_{k}})\right).
\end{array}
$$
Further on,
$$
\phi(e_{s_{2}}\times\tilde e_{t_{1}}(\hat\theta_{nm}))=
e_{s_{2}}\times\tilde e_{t_{1}}(\phi(\hat\theta_{nm}))=
e_{s_{2}}\times\tilde e_{t_{1}}(\sum v^*_L\otimes v_L),
$$
where $L$ runs over all the sequences  of length $nm$
composed from the integers 1 to $n$. But, as is not difficult to
see,
$$
e_{s_{2}}\times\tilde e_{t_{1}}(\sum v^*_L\otimes v_L)=
const\cdot e_{s_{2}}(v^*_{I_{m}})\otimes\tilde
e_{t_{1}}(v_{I_{m}}),
$$
hence,
$$
X_-X_+e_{s_{2}}(v^*_{I_{m}})\otimes\tilde
e_{t_{1}}(v_{I_{m}})=const\cdot e_{s_{2}}\times e_{s_{1}}
(\hat\theta_{nm}). 
$$
Therefore,
$$
w=const\cdot \mathop{\sum}\limits_{\sigma, \tau}\eps(\tau)\sigma\times
\tau\left (e_{s_{1}}(v^*_{I_{k}})\otimes
e_{s_{2}}\times\tilde e_{t_{1}}(\hat\theta_{nm})\otimes
e_{t_{2}}(v_{J})\right)=
e_{s}\times\tilde e_{t}(v^*_{I_{k}}\otimes\hat\theta_{nm}
\otimes v_{J})
$$
which proves the theorem. \end{proof}

Proof of the following theorem is similar.

\ssbegin{3.4}{Theorem} The element $e_{s}\times\tilde 
e_{t}(v_{I_{k}}\otimes\theta_{nm}\otimes v_{J_{k}}^*)$ in 
$T^{nm+kn, nm+km}(V)$ is $\fsl(V)$-invariant.
\end{Theorem}

\ssbegin{3.5}{Corollary} Let $L$ eb the sequence with elements from 
$\Mee$. Set 
$$
p(L)=\sum p(l_{i}),\quad \alpha(L, 
L)=\mathop{\sum}\limits_{i<j}p(l_{i})p(l_{j}). 
$$
Under notations of Theorems $3.3$, $3.4$ the
invariant elements can be expressed in the form
$$
e_{s}\times\tilde
e_{t}(v^*_{I_{k}}\otimes\theta_{nm}\otimes v_{J_{k}})=
\mathop{\sum}\limits_L(-1)^{p(L)+\alpha(L, L)}e_{s}(v^*_{I_{k}}\otimes v^*_L)\otimes\tilde
e_{t}(v_L\otimes v_{J_{k}})\eqno{(3.5.1)}
$$
and 
$$
e_{s}\times\tilde
e_{t}(v_{I_{k}}\otimes\theta_{nm}^*\otimes v^*_{J_{k}})=
\mathop{\sum}\limits_L(-1)^{\alpha(L, L)}e_{s}(v_{I_{k}}\otimes v_L)\otimes\tilde
e_{t}(v^*_L\otimes v^*_{J_{k}}),\eqno{(3.5.2)}
$$
where the sums run over all the sequences $L$ of length nm with
elements from $R_V$.
\end{Corollary}

\begin{proof} It is easy to verify that $\hat \theta_{nm}=\mathop{\sum}\limits_L 
(-1)^{\alpha(L, L)+p(L)}v^*_L\otimes v_L$, which immediately implies 
(3.5.1).  Formula (3.5.2)  is similarly proved. \end{proof}

\ssec{3.6} Recall the definition of $R_V$, $R_U$ and $R_W$ (0.4) and 
(0.5).  For any sequences $I$ and $J$ denote by $I*J$ the sequence 
obtained by ascribing $J$ at the end of $I$.  Let now $I$ be the 
sequence of length $(k+m)n$ with elements from $R_U$, let $J$ be the 
sequence of length $(k+n)m$ with elements from $R_W$, let $\hat I$ be 
the sequence of length $(k+n)m$ with elements from $R_U$ and $\hat J$ 
be the sequence of length $(k+m)n$ with elements from $R_W$.

For any sequence $L$ of length $nm$ with
elements from $R_V$ we define:
$$
\renewcommand{\arraystretch}{1.4}
\begin{array}{l}
\tilde P_s(I, I_k*L)\in S^{\bcdot}(U\otimes V),\quad
\tilde P_t(L*J_k, J)\in S^{\bcdot}(V^*\otimes W);\\
P_t(\hat I, L*J_k)\in S^{\bcdot}(U\otimes V), \quad P_s(I_k*L, \hat J)\in
S^{\bcdot}(V^*\otimes W).
\end{array}
$$

\begin{Theorem} The algebra of $\fsl(V)$-invariant elements in 
$\fA^{p, q}_{k, l}$ is generated by the elements

{\em i)} $(v^*_r, v_s)$, where $r\in R_U$ and $s\in
R_W$;

{\em ii)} $F_k(I, J)=\mathop{\sum}\limits_L(-1)^{\alpha(L, L)}\tilde 
P_s(I, I_k*L)\tilde P_t(L*J_k, J)$, where $I$ is an $s$-semistandard 
sequence and $J$ is a $t$-semistandard one;

{\em iii)} $F_{-k}(\hat I, \hat J)=
\mathop{\sum}\limits_L(-1)^{\alpha(L, L)+p(L)(p(\hat I)+p(\hat J))}
P_s(I_k*L, \hat J) P_t(\hat I, L*J_k)$,

\noindent where $\hat I$ is an $s$-semistandard 
sequence, $\hat J$ is a $t$-semistandard one and $L$ runs over all 
the sequences of length $nm$ with
elements from $R_V$.
\end{Theorem}

\begin{proof} For Young tableaux $\lambda$ and $\mu$ we
have
$$
(V^\lambda\otimes
V^*{}^\mu)^{\fsl(V)}=\Hom_{\fsl(V)}(V^\mu, V^\lambda).
$$
The dimension of this space is equal to either 0 or 1. It is equal
to 1 only if $\lambda=\mu$ or both of them contain a
$n\times m$ rectangle and $\lambda_i=\mu_i+k$ for $i=
1, \dots , n$ and $\lambda'_j=\mu'_j+k$ for $j=1, \dots , m$
and any $k\in \Zee$.

To prove the theorem, it suffices to show that for the above $\lambda$ 
and $\mu$ the module $V^\lambda\otimes V^*{}^\mu$ containes an 
invariant which can be expressed via the invariants listed in the 
theorem.  By \cite{S2}  such an invariant exists.  Under the canonical 
homomorphism of the tensor algebra onto the symmetric one, the 
invariants of the form i)--iii)  turn into a system of generators.  
Theorem is proved.  \end{proof}

\ssec{3.7} To the invariant element in $T^{n(m+k), m(n+k)}(V)$ there 
corresponds an invariant operator $T^{m(n+k)}(V)\tto T^{n(m+k)}(V)$.  
To describe it, observe that $C_t$ can be represented as 
$C_t=\mathop{\coprod}\limits_{\pi\in Z}(C_{t_{1}}\times C_{t_{2}})\pi$, 
the decomposition into right cosets relative  
the product of the column stabilizers of tableaux $t_{1}$ and 
$t_{2}$; let $Z$ be a collection of their representatives.  
Define
$$
D_{J_{k}}: T^{m(n+k)}(V)\tto T^{nm}(V),\quad
D_{J_{k}}(v_1\otimes v_2)=(-1)^{p(J_{k})p(v_{1})} v_1\cdot
v^*_{J_{k}}(v_2).
$$

\begin{Lemma} Let $\cL$ be an invariant operator
corresponding to $e_{s}\times\tilde
e_{t}(v_{I_{k}}\otimes\theta_{nm}\otimes v^*_{J_{k}})$.
Then
$$
\cL(e_t(v_L))=const\cdot
\cL(v_L)=e_{s}(v_{I_{k}}\otimes
D^*{}_{J_{k}}e_{t_{2}}\mathop{\sum}\limits_{\pi\in
Z}\eps(\pi)\pi v_L).\eqno{(3.7)}
$$
\end{Lemma}

\begin{proof} To $\theta_{nm}$ there corresponds the identity operator 
$\id: V^{\otimes nm}\tto V^{\otimes nm}$; hence, to 
$\theta_{nm}\otimes v^*_{J_{k}}$ there corresponds the operator 
$D_{J_{k}}: V^{\otimes m(n+k)}\tto V^{\otimes nm}$ and to 
$v_{I_{k}}\otimes\theta_{nm}\otimes v^*_{J_{k}}$ there corresponds the 
operator $v_{I_{k}}\otimes D_{J_{k}}$; finally, to 
$e_{s}\times\tilde e_{t}(v_{I_{k}}\otimes\theta_{nm}\otimes 
v^*_{J_{k}})$ there corresponds the operator $e_{s}(v_{I_{k}}\otimes 
D_{J_{k}})e_{t}$.  Hence,
$$
\renewcommand{\arraystretch}{1.4}
\begin{array}{l}
    \cL(e_{t}(v_L))=e_{s}v_{I_{k}}\otimes D_{J_{k}})e_{t}^ 2(v_L=
c_{1}e_{s}v_{I_{k}}\otimes D_{J_{k}})e_{t}(v_L)=
c_{1}\cdot \cL(v_L)=\\
c_{1}e_{s}(v_{I_{k}}\otimes D_{J_{k}}\mathop{\sum}\limits_\pi 
e_{t_{1}}e_{t_{2}}\eps(\pi)\pi v_L)= c_{1}e_{s}(v_{I_{k}}\otimes 
e_{t_{1}}D_{J_{k}}e_{t_{2}}\mathop{\sum}\limits_\pi \eps(\pi)\pi 
v_L)=\\
c_{1}e_{s}e_{t_{1}}(v_{I_{k}}\otimes 
D_{J_{k}}e_{t_{2}}\mathop{\sum}\limits_\pi \eps(\pi)\pi v_L)= 
c_{1}c_{2}e_{s}(v_{I_{k}}\otimes 
D_{J_{k}}e_{t_{2}}\mathop{\sum}\limits_\pi \eps(\pi)\pi v_L).
\end{array}
$$
The last equality follows from 
$e_{s}e_{t_{1}}=c_{2}e_{s}$. \end{proof}

\ssec{3.8} Let us consider the case $k=1$ in more detail.  Let $L$ be a sequence 
of length $nm+m$ with elements from $R_V$, considered as a $t$-tableau.

In each column $L$, mark an ``odd" element so that all the elements 
marked, say, $l=(l_1, \dots , l_m)$, are distinct.  The pair $(L, l)$ 
will be called a {\it marked} tableau.  Introduce the following 
notations: $c_i$ for the parity of the $i$-th column, $d_i$ for the 
parity of the last elelment in the $i$-th column, $b_i$ for the parity 
of the column under the $i$-th marked element, $|b_i|$ for the number 
of elements in the $i$-th column under the $i$-th marked element, 
$\eps(l)$ for the sign of the permutation $l=(l_1, \dots , l_m)$ and 
$\eps(L, l)=(-1)^{q(L)}\eps(l)$; set further
$$
\eps(L)=c_2+c_4+\dots +d_2+d_4+\dots; \; \; 
q(L)=b_1+|b_1|+b_2+|b_2|+\dots .
$$

\begin{Theorem} The invariant operator is of the form
$$
\cL(e_{t}(v_L))=const\cdot
\cL(v_L)=const\cdot\eps(L)\mathop{\sum}\limits_{(L, l)}
\eps(L, l)e_{s}(v_{I_{1}}\otimes v_{L\setminus l}),\eqno{(3.8)}
$$
where the constant factor does not depend on $L$.\end{Theorem}

\begin{proof} Since for the representatives of the cosets of 
$\fS_{n+1}/\fS_{n}$ we can take a collection of 
cycles, we may assume in formula (3.7)  that
$$
\pi=\pi_1\dots\pi_m,\;\; \pi_i\pi_j=\pi_j\pi_i\text{ for any }i,
j.
$$
Hence, $\pi^ 2=1$.  Furhter on, $D_{J_{1}} e_{t_{2}}\sum\eps(\pi)\pi 
v_{L}\neq 0$ if and only if the last row of $\pi L$ for some $\pi$ is, 
up to a permutation, a permutation of $\{\bar 1, \dots , 
\bar m\}$.

The set of marked tableau $(L, l)$ is in one-to-one correspondence 
with the set of pairs $(L, \pi)$ such that the last row of $\pi L$ is, 
up to a permutation, $\{\bar 1, \dots , \bar m\}$.  Indeed, from the 
pair $(L, l)$ determine $\pi=\pi_1\dots\pi_m$, where $\pi_i$ is the 
cycle that shifts the elements under the $i$-th marked one one cell up 
along the column and places the marked one at the bottom.  If the 
marked element lies in the last row, we set $\pi_i=1$.

And, the other way round, given $\pi$, we mark $\pi(k_1)$, \dots 
, $\pi(k_m)$, where $(k_1, \dots , k_m)$ is the last row of $L$.  
Hence, (3.7)  implies that
$$
\cL(v_L)=\mathop{\sum}\limits_{(L, l)}
\delta(L, l)e_{s}(v_{I_{1}}\otimes v_{L\setminus l}),
$$
where $\delta(L, l)$ is a sign depending on $(L, l)$.  Direct 
calculations of this sign lead us to (3.8).  \end{proof}

\section*{\S 4. The absolute invariants of \protect $\fosp(V)$}

Let $A=U(\fosp(V))[\eps]$ be the central extension with the
only extra relation $\eps ^2=1$.

On $A$, introduce the coalgebra structure making use of that on 
$U(\fosp(V))$ and setting $\eps\mapsto \eps\otimes \eps$.  Assuming 
that $\eps$ acts on $V$ as the scalar operator of multiplication by 
$-1$, we may consider $V$ as an $A$-module.  Using the coalgebra 
structure on $A$, one can determine a natural $A$-action in $T^{p, 
q}(V)$ and $\fA^{p, q}_{k, l}$.  So, we can speak about $A$-invariants 
in these modules.

\ssbegin{4.1}{Lemma} Let $\fgl(V)=\fg=\fg_-\oplus \fg_0\oplus\fg_+$, as 
in \S 3  and let $M$ be a $\fg_{0}$-module.  Set $\fg_+M=0$.  There is an 
isomorphism of $\fosp(V)$-modules
$$
\ind^{\fgl(V)}_{\fg_0\oplus \fg_+}(M)\simeq 
\ind^{\fosp(V)}_{\fosp(V)_\ev}(M).\eqno{(4.1)}
$$
\end{Lemma}

\begin{proof} Cf.  \cite{S2}, Lemma 5.1.   
\end{proof}

\ssbegin{4.2}{Lemma} Let $\fg$ be a Lie superalgebra and the
representation of $\fg_\ev$ in the maximal exterior power of 
$\fg_\od$
is trivial. Then there is an isomorphism of vector spaces
$$ 
\ind^{\fg}_{\fg_\ev}(M)^{\fg}\simeq M^{\fg_\ev}.\eqno{(4.2)}
$$
\end{Lemma}

\begin{proof} Cf.  \cite{S2}, Lemma 5.2. 
\end{proof}

\begin{Remark} Statements similar to Lemmas 4.1 and 4.2 hold also for 
$U(\fosp(V))[\eps]$-modules.  One can refine Lemma 4.2 and prove that 
if $v_0\in M$ is $\fg_\ev$-invariant, then the corresponding to it 
$\fg$-invariant vector is of the form $\xi_1\dots\xi_nv_0+$ terms of 
lesser degree.
\end{Remark}

\ssec{4.3} The presence of an even $\fosp(V)$- and $A$-invariant form 
on $V$ determins an isomorphism of $A$-modules and algebras $\fA^{p, 
q}_{k, l}=\fA^{p+k, q+l}$.  Therefore, we may assume that $k=l=0$.  By 
definition, the Lie superalgebra $\fosp(V)$ preserves the vector
$$
\mathop{\sum}\limits_{i=1}^ne_i^*\otimes 
e_{n-i+1}^*+\mathop{\sum}\limits_{j=1}^r(e_{\overline{m-j+1}}^*\otimes 
e_{\bar j}^*-e_{\bar j}^*\otimes e_{\overline{m-j+1}}^*),
$$ 
where $\dim V=(n|2r)$. Therefore, the scalar products
$$
(v_s,v_t)=\mathop{\sum}\limits_{i=1}^nx_{is}^*x_{n-i+1, 
t}^*+(-1)^{p(s)}\mathop{\sum}\limits_{j=1}^r(x_{\overline{m-j+1}, 
s}^*x_{\bar j, t}^*-x_{\bar j, s}^*x_{\overline{m-j+1, t}}^*), 
\eqno{(4.3)}
$$
where $s, t\in R_W$, are $\fosp(V)$- and $A$-invariants.

\begin{Theorem} The algebra of $A$-invariant elements
in $\fA^{p, q}=S^{\bcdot}(V^*\otimes W)$ is generated by the elements $(v_s,
v_t)$ for $s, t\in R_W$. \end{Theorem}

\begin{proof} Cf.  \cite{S2}, Theorem 5.3. 
\end{proof}

\ssec{4.4} Let $I$ be a sequence of length $2k$ with elements from $R_W$.
Determine an element $X(I)\in S^{\bcdot}(S^2(W))$ by setting
$$
X(I)=x_{i_{1}i_{2}}\dots x_{i_{2k-1}i_{2k}},
$$
where $x_{ij}$ is the canonical image of the element
$w_i\otimes w_j\in S^2(W)$. 

Let $t$ be a tableau of order $2k$ with rows of even lengths.
Then an ``even Pfaffian" is defined:
$$
Pf_t(I)=\mathop{\sum}\limits_{\tau\in C_{t}, \; \sigma \in 
R_{t}}\eps(\tau)c(I, (\sigma \tau)^{-1})X(\sigma \tau I).\eqno{(4.4)}
$$

\begin{Theorem} {\em a)} $S^{\bcdot}(S^2(W))=\oplus W^{\lambda}$, 
where the length of each row of $\lambda$ is even.

{\em b)} Let $t$ be a $\lambda$-tableau filled in along rows with the 
numbers $1$, $2$, \dots .  Then the family $Pf_t(I)$ for the 
$t$-standard sequences $I$ is a basis of $W^{\lambda}$.
\end{Theorem}

\begin{proof} On $T^{2k}(W)=W^{\otimes 2k}$ the group $\fS_{2k}$ and 
its subgroup $G_k=\fS_{k}\circ \Zee_2^ k$ naturally act; namely, 
$\fS_{k}$ permutes pairs $(2i-1, 2i)$ whereas $\Zee_2^ k$ permutes 
inside each pair.

Clearly, $S^k(S^2(W))=T^{2k}(W)^{G_{k}}$. 

But, on the other hand, $T^{2k}(W)=\oplus\; 
S^\lambda\otimes W^\lambda$, so $T^{2k}(W)^{G_{k}}=\oplus\; 
(S^\lambda)^{G_{k}}\otimes W^\lambda$. Hence, in the decomposition
of $S^k(S^2(W))$ only enter $W^\lambda$ for which
$(S^\lambda)^{G_{k}}\neq 0$ and their multiplicity is equal to
$\dim (S^\lambda)^{G_{k}}$.

But
$$
(S^\lambda)^{G_{k}}=\Hom_{G_{k}}(\ind^{\fS_{k}}_{G_{k}}(\id),
S^\lambda),
$$
so the multiplicity of  $W^\lambda$ in $S^k(S^2(W))$ is equal
to that of  $S^\lambda$ in $\ind^{\fS_{k}}_{G_{k}}(\id)$. By
\cite{H} it is equal to 1 if the lengths of all rows of $\lambda$ are even
and 0 otherwise. This proves a).

b) Consider now the natural map $T^{2k}(W)\tto S^k(S^2(W))$.  For the 
tableau $t$ from the conditions of the theorem and the sequence $I$ 
the vectors $e_t(w_I)$ form a basis of $W^\lambda$.  So the images of 
these vectors (which are exactly the $Pf_t(I)$) form a basis of 
$W^\lambda\subset S^k(S^2(W))$.
\end{proof}

\ssec{4.5} Consider the algebra homomorphism
$$
S^{\bcdot}(S^2*W))\tto S^{\bcdot}(V^*\otimes W),\quad x_{st}\mapsto (v_s, 
v_t).\eqno{(4.5)}
$$
Its kernel is the ideal of relations between scalar products.

\begin{Theorem} The ideal of relations between scalar
products is generated by polynomials $Pf_t(I)$, where $t$ is a
$(2r+2)\times (n+1)$ rectangle filled in along rows and $I$ is a
$t$-standard sequence with elements from $R_W$. \end{Theorem}

\begin{proof} By Theorem 4.3
$$
S(V^*\otimes W)^A=
\mathop{\oplus}\limits_\lambda(V^*{}^\lambda)^A\otimes W^\lambda=
\mathop{\oplus}\limits_{\lambda_{n+1}\leq 2r}W^\lambda.
$$
Hence, the kernel of homomorphism (4.5) is equal to 
$\mathop{\oplus}\limits_{\lambda_{n+1}\geq 2r+2}W^\lambda$.  Let us 
prove that it is contained in the ideal generated by $W^\lambda$, 
where $\lambda$ is a $(2r+2)\times (n+1)$ rectangle.

Let $\mu\supset\lambda$ and $e_s$ the corresponding
idempotent. Then $e_s=\sum\tau_ie_t\sigma_j$. Hence,
$$
e_s(J)=\mathop{\sum}\limits_{}\tau_ie_t(\sigma_jJ)= 
\mathop{\sum}\limits_{}e_{\tau_{i}t}(\tau_{i}\sigma_jJ).
$$
Thus, $Pf_s(J)=\mathop{\sum}\limits_{i,j}f_{ij}Pf_{\tau_{i}t}(J_{ij})$ 
and we are done.  \end{proof}

\section*{\S 5. The relative invariants of \protect $\fosp(V)$}

The invariants of $\fosp(V)$ are, first of all, the ones
generated by scalar products. To describe the other invariants,
let us describe a certain invariant in the tensor algebra. Let
$\dim V= n|m$. For $i\in R_V$
define $\tilde i$ by setting
$$
\tilde i= n-i+1 \text{ if $i$ is ``even" and $\overline{m-i+1}$ if $i$ is
``odd"}.  
$$
Let $I=i_1i_2\dots i_{2p}$ be a sequence of even length with elements 
from $R_V$ and $I^*$ the set consisting of the pairs $(i_{2\alpha-1}, 
i_{2\alpha})$ for $\alpha\leq p$ such that $\tilde i_{2\alpha-1}\neq 
i_{2\alpha}$.  Let $t$ be a rectangular $n\times m$ tableau 
consequtively filled in along columns from left to right and $I$ a 
sequence with elements from $R_V$.  Let us fill in the tableau $t$ 
with elements from $I$: replace $\alpha$ with $i_\alpha$.  Let $\cT$ 
be the set of sequences $I$ such that all the rows of $t$ except the 
last row are of the form
$$
i_1\tilde i_1\dots i_{r}\tilde i_r\text{ for }r=\frac 12 m,
$$
while the last row $J$ should be such that if $j\in\hat J$, then 
$\tilde j\in\hat J$ and $\hat J$ consist of pairwise
distinct ``odd'' elements.

Let $I\in\cT$. Set $r=\frac12m$, let $\nu$ be the
total amount of marked pairs from the last row consisting of
pairwise conjugate ``odd" elements that do not belong to
$N(L)$; let $n_1$, \dots , $n_\nu$ the multiplicities with which
these pairs enter the last row and $N=n_1+ \dots+n_\nu$; let
$\sigma_l$ be the $l$-th elementary symmetric function. Set
$$ 
K(I)=\mathop{\sum}\limits_{q=s}^{s+\nu}(N+1)^r2^{r-q}(r-q)!  
N^q\sigma_{q-s}(n_1, \dots , n_\nu)
$$
and $d(I)=d(I_1)d(I_3)\dots d(I_{2r-1})$,  where $d(J)=(-1)^{\alpha(J, J)}$, 
see $(1.2)$??.

\ssbegin{5.1}{Theorem} In $V^{\otimes  n(m+1)}$ lies an
$\fosp(V)$-invariant element
$$
\nabla_{m+1}=\mathop{\sum}\limits_{I\in 
\cT}d(I)K(I)e_s(v_{I_{1}}\otimes v_{I}).\eqno{(5.1.1)}
$$
\end{Theorem}

\begin{proof} Set
$$
c(i, \tilde i)=\cases 1&\text{ if $p(i)=0$ or
$i<\tilde i$ and  $p(i)=1$}\cr
- 1&\text{ if $i>\tilde i$ and  $p(i)=1$}.\endcases
$$
The map 
$$
V\tto V^*, \quad e_i\mapsto c(i, \tilde i)e^*_{\tilde i}
$$ 
is an isomorphism induced by the invariant bilinear form and $\tilde 
\theta_2=\mathop{\sum}\limits_{i\in R_V} c(i, \tilde i)e_i\otimes 
e_{\tilde i}$ is an $\fosp(V)$-invariant.

Let $t$ be a rectangular $(n+1)\times m$ tableau as in Theorem 3.8 and 
$J$ a $t$-sequence such that after being filled each row $J$ is of the 
form $j_1\tilde j_1\dots j_r\tilde j_r$.  Denote by $\cT_1$ the set of 
such sequences $J$.  Then
$$
\tilde \theta=\theta_2^{\otimes \frac12
(n+1)m}=\mathop{\sum}\limits_{J\in\cT_1} d(J)c(J)v_J,\eqno{(5.1.2)}
$$
where $J_1$, \dots , $J_{2r-1}$ are the columns of the tableau
$t$ and where
$$
d(J)=d(J_1)d(J_3)\dots d(J_{2r-1})\text{ while }c(J)=c(J_1)c(J_3)\dots
c(J_{2r-1})
$$ 
whereas $c(J_\alpha)=\mathop{\prod}\limits_{i\in J_{\alpha}} c(i, 
\tilde i)$.
 
 The element (5.1.2) is an
$\fosp(V)$-invariant; having applied to it the operator $\cL$
from Theorem 3.8 we get another $\fosp(V)$-invariant:
$$
\cL(\tilde \theta)=\mathop{\sum}\limits_{J} d(J)c(J)\cL(v_J)=
\mathop{\sum}\limits_{J, l} d(J)c(J)\eps(J)\eps(J, l)e_s(v_{I_{1}}\otimes
v_{J\setminus l}),\eqno{(5.1.3)} 
$$
where 
$$
\eps (J)=\mathop{\prod}\limits_{1\leq i\leq r}\eps(J_{2i-1}*J_{2i}),\quad 
\eps (J, l)=\sign (l)\mathop{\prod}\limits_{1\leq i\leq
r}\eps(J_{2i-1}*J_{2i}, l_{2i-1}*l_{2i}).
$$
For the collection $(J_{2i-1}, J_{2i}, l_{2i-1}, l_{2i})$ define the
sequence $(I_{2i-1}, I_{2i})$ as follows: if $ l_{2i-1}$
and $l_{2i}$ lie in the same row just strike them out,  if $ l_{2i-1}$
and $l_{2i}$ lie in distinct rows we strike them out and
place their conjugates, $\tilde l_{2i-1}$
and $\tilde l_{2i}$, in the last row in the same columns. The sequence $I$
takes the form $(I_1, I_2, \dots ,  I_{2r-1}, I_{2r})$. It is not
difficult to verify that
$$
e_s(v_{I_{1}}\otimes
v_{J\setminus l})=\sign(l)\eps(J, l)e_s(v_{I_{1}}\otimes
v_{I})\text{ and }d(J)=(-1)^rd(I)\eps(J). 
$$
Therefore, 
$$
\cL(\tilde \theta)=(-1)^r\mathop{\sum}\limits_{J, 
l}\sign(l)c(J)d(I)e_s(v_{I_{1}}\otimes v_{I}).
$$
The constant factor, the sign, can be, clearly, replaced with a 1.  If 
$I$ is of the above form, then in the last row for some values of $i$ 
the pairs $(l_{2i-1}, l_{2i})$ are conjugate whereas all the remaining 
values of $i$ are odd and pairwise distinct, call them $I^*=\{k_1, \dots , 
k_{2p}\}$.  Set
$$
\hat c(I)=\sign(k_1, \dots , 
k_{2p})(-1)^{p}\mathop{\prod}\limits_{c(i_{2\alpha -1}, i_{2\alpha 
}),\neq 0}c(i_{2\alpha -1}, i_{2\alpha}),
$$
where $\sign(k_1, \dots , 
k_{2p})$ is the sign of the permutation. Then $c(J)\eps(l)=\hat c(I)$ 
and, therefore,
$$
\cL(\theta_{m+1})=\mathop{\sum}\limits_{J, 
l}c(I)d(I)e_s(v_{I_{1}}\otimes v_{I}).
$$
To complete the proof, it suffices to calculate the number of
pairs$(J, l)$ that give the sequence $I$ which leads to formula
(5.1.1). \end{proof}

\ssbegin{5.2}{Theorem} The algebra of $\fosp(V)$-invariants is
generated by the polynomials

{\em i)} $(v_s, v_t)$ for $s, t\in R_W$ and 

{\em ii)} $R(J)=\mathop{\sum}\limits_{I}d(I)K(I)Pf_s(I_1*I, J)$ for 
every $I\in \cT$ and every $s$-standard sequence $J$ with elements 
from $R_W$.
\end{Theorem}

\begin{proof} Let $f$ be an $\fosp(V)$-invariant which is not
$A$-invariant. Let $f$ depend on $n-1$ even and $2r$ odd
generic vectors $v_1$, \dots , $v_{\overline{2r}}$. Then there
exists a $g\in\OSp(V\otimes A)$ such that $g\Span(v_1, \dots ,
v_{\overline{2r}})=\Span (e_1, \dots , e_{\overline{2r}})$. 

Let $he_n=-e_n$ whereas $he_i=e_i$ for $i\neq n$. Then $\ber
(h)=-1$ and $f(hg\cL)=-f(g\cL)$. But, on the other hand, 
$f(hg\cL)=f(g\cL)$, hence, $f=0$. This means that
$\fosp(V)$-invariants other than scalar products may only be of
type $\lambda$ corresponding to a typical module. So, in the
same vein as for $A$-invariants, we see that $\dim
(V^*{}^\lambda)^{\fosp(V)}=1$ if $\lambda$ is typical and its
first $n$ rows are of odd lengths whereas the remaining
rows are of even lengths. If we do not consider the scalar
products, then for the other (atypical) $\lambda$ there are
no invariants. 

Under the canonical homomorphism
$T^k(V^*)\otimes T^k(W)\tto S^k(V^*\otimes W)$ the module 
$V^*{}^\lambda\otimes W^\lambda$ turns into its copy and a
basis of the first copy becomes a basis of the second one. This
shows that if  $\lambda$ is an $n\times (2r+1)$ rectangle,
then the polynomials $R(J)$ from the theorem constitute a
basis of $V^*{}^\lambda\otimes W^\lambda$, a subspace of
$S^k(V^*{}\otimes W)$. 

For an arbitrary $\lambda$ containing
an $n\times (2r+1)$ rectangle we apply the same arguments
as in the proof of Theorem 2.2.
\end{proof}

\section*{\S 6. The invariants of \protect $\fpe(V)$}

Suppose $\dim V=(n|n)$, the $e_i^*$ is a basis of
$V_\ev$ and  the $e_{\bar i}^*$ be the dual basis of
$V_\od$ with respect to an odd nondegenerate form on
$V$. Then $\fpe(V)$ preserves the tensor $\sum
(e_i^*\otimes e_{\bar i}^*+e_{\bar i}^*\otimes e_i^*)$.

Observe that the scalar products 
$$
(v_s,  v_t)=\mathop{\sum}\limits_{}(-1)^{p(s)}
(x_{is}^*\otimes e_{\bar i, t}^*+e_{\bar i, s}^*\otimes e_{it}^*)
\text{ for any } s, t\in S
$$
are $\fpe(V)$-invariants. Moreover, the presence of the odd
form determines an isomorphism of algebras and
$\fpe(V)$-modules $\fA^{p, q}_{k, l}=\fA^{p+l, q+k}$, so, as for the
orthosymplectic case, we may assume that $k=l=0$.

The compatible $\Zee$-grading of $\fgl(V)$ induces compatible
$\Zee$-gradings of $\fpe(V)$ and $\fspe(V)$:
$$
\fg=\fg_-\oplus\fg_0\oplus \fg_+,\text{ where
$\fg_-=\Lambda^2(V)$, $\fg_+=S^2(V^*)$ and $\fg_0=\fgl(V)$
or $\fsl(V)$}
$$
(There is also another, isomorphic, representation which we
will not use in this paper: 
$$
\fg=\fg_-\oplus\fg_0\oplus \fg_+,\text{ where
$\fg_-=\Lambda^2(V^*)$, $\fg_+=S^2(V)$ and $\fg_0=\fgl(V)$
or $\fsl(V)$}.)
$$
Let $X_\alpha$, $1\leq \alpha\leq \frac12n(n+1)$, be a basis of
$\fg_+$ and let $Y_\beta$,
$1\leq \beta\leq \frac12n(n-1)$, be a basis of $\fg_-$. Set
$$
X_+=\mathop{\prod}\limits_{1\leq \alpha\leq \frac12n(n+1)}X_\alpha, \quad 
Y_-=\mathop{\prod}\limits_{1\leq \beta\leq \frac12n(n-1)}Y_\beta.
$$

Observe that the weight of $X^+$ with respect to the Cartan subalgebra 
is equal to $(n+1)\sum\eps_i$ and the weight of $Y^-$ is equal to 
$-(n-1)\sum\eps_i$.

\ssbegin{6.1}{Lemma} Let $L=\ind^\fg_{\fg_0\oplus\fg_+}(M)=
\ind^\fg_{\fg_0\oplus\fg_-}(N)$ be a typical irreducible
$\fg=\fgl(V)$-module. There exists an isomorphism of vector
spaces $$ L^{\fspe(V)}=M^{\fspe(V)_{\ev}}=N^{\fspe(V)_{\ev}}
$$
given by the formulas
$$
M\tto L, \; m\mapsto Y^-m\text{ and }N\tto L, \; n\mapsto
X^+n
$$
the inverse map being given by the formulas
$$
L\tto M, \; l\mapsto X^-l\text{ and }L\tto N, \; l\mapsto
Y^+l.
$$
\end{Lemma}

\begin{proof} Consider the two gradings of $L$:
$$
L^+_k=\Span(f(X_\alpha)n: n\in N \text{ and } \deg f=k)
$$
 and 
 $$
L^-_k=\Span(f(Y_\beta)m: m\in M\text{ and } \deg f=k).
$$
It is clear that $L^+_k=L^-_{n^{2}-k}$.

If $l$ is a $\fspe(V)$-invariant, then $X_\alpha l=0$ (for 
$1\leq \alpha\leq \frac12n(n+1)$)  and $l=X^+f(X_\alpha)n$ for $n\in N$.  Therefore, 
$l=\mathop{\sum}\limits_{r\geq \frac12n(n+1)}l_r^+$, where $l_r^+\in 
L_r^+$.  We similarly establish that $l=\mathop{\sum}\limits_{1\leq s\leq 
\frac12n(n-1)}l_s^-$, where $l_s^-\in L_s^-$.  Hence, 
$\mathop{\sum}\limits_{r\geq 
\frac12n(n+1)}l_r^+=\mathop{\sum}\limits_{1\leq s\leq \frac12n(n-1)}l_s^-$.  
Taking into account the equality $L^+_k=L^-_{n^{2}-k}$ we deduce that 
$l\in L^+_{\frac12n(n+1)}=L^-_{\frac12n(n-1)}$ and $l=X^+n=Y^-m$ for 
some $n\in N$ and $m\in M$.  Moreover, it is clear that $m$ and $n$ 
are $\fspe(V)_\ev=\fsl(V_\ev)$-invariants.

Conversely, if $m$ and $n$ are $\fsl(V_\ev)$-invariants, then a
direct check shows that $X^+n$ and $Y^-m$ are
$\fspe(V)$-invariants. \end{proof}

\ssbegin{6.2}{Theorem} The algebra of $\fpe(V)$-invariants is
generated by the scalar products $(v_s, v_t)$ for $s, t\in R_W$.
\end{Theorem}

\begin{proof} See\cite{S2}, sec. 6.2. \end{proof}

Let $I$ be a sequence of length $2k$ composed of elements
from $R_W$. Determine the element $Y(I)\in
E^{\bcdot}(S^2(W))=S^{\bcdot}(\Pi(S^2(W)))$ by setting
$$
Y(I)=(-1)^{\beta}y_{i_{1}i_{2}}\dots y_{i_{2k-1}i_{2k}},
$$
where $y_{ij}$ is the canonical image of the element $\omega_i\otimes 
\omega _j$ and $\beta=\mathop{\sum}\limits_{1\leq \alpha\leq 
k}(k-\alpha)(i_{2\alpha-1}+i_{2\alpha})$.

\ssec{6.3} Let $\lambda$ be a partition of the form $(\alpha_{1},
\dots , \alpha_{p}, \alpha_{1}-1, \dots, \alpha_{p}-1)$ in
Frobenius' notations (see \cite{M}) and $t$ be a tableau of the form
$\lambda$ filled in so that the underdiagonal columns the
diagonal cells including are filled in consequtively with odd
(ili ``odd"?) numbers while the rows to the right of the diagonal
are consequtively occupied by even numbers. For a tableau of
such a form and a sequence $I$ the ``periplectic" Pfaffian is
defined: 
$$
PPf_t(I)=\mathop{\sum}\limits_{\tau\in C_t, \;\sigma\in 
R_t}\eps(\tau)c(I, (\sigma\tau)^{-1})Y(\sigma\tau I).\eqno{(6.3.1)}
$$

\ssbegin{6.3.1}{Theorem} For the above tableau $t$ the family
$PPf_t(I)$ for the $t$-standard sequences $I$ is a basis in the
module $W^\lambda\subset E^{\bcdot} (S^2(W))$.
\end{Theorem}

\begin{proof} From the theory of $\lambda$-rings it follows
that $E^{\bcdot} (S^2(W))=\oplus W^\lambda$, where the sum runs over
the $\lambda$ of the above described form. One can easily
verify that for the tableau as indicated in the formulation of
the theorem and a $t$-standard sequence $I$ the image
$e_t(w_I)$ in $E^{\bcdot} (S^2(W))$ is nonzero. Hence, for a fixed
tableau $t$ the canonical map $T^{2k}(W)\tto E^k(S^2(W))$
performes an isomorphism of $e_t(T^{2k}(W))$ with
$W^\lambda\subset E^k(S^2(W))$. This implies the theorem.
\end{proof}

\ssec{6.3.2} Consider now an algebra homomorphism 
$$
E^{\bcdot} (S^2(W))\tto S^{\bcdot} (V^*\otimes  W),\quad y_{st}\mapsto
(v_s, v_t). \eqno{(6.3.2)} 
$$

\begin{Theorem} The kernel of $(6.3.2)$ is generated by 
polynomials $PPf_t(I)$, where $t$ is of the form of a 
$(n+1)\times(n+2)$ rectangle and is filled in as described in the 
previous section and $I$ is a $t$-standard sequence with elements from 
$R_W$.
\end{Theorem}
    
\begin{proof} Clearly,
$$
(S^k(V^*\otimes W))^{\fpe(V)}= (\mathop{\oplus}\limits_{\lambda: 
\lambda _{n+1}\leq n}V^*{}^{\lambda}\otimes 
W^*{}^{\lambda})^{\fpe(V)}= \mathop{\oplus}\limits_{\lambda: \lambda 
_{n+1}\leq n}W^*{}^{\lambda},
$$
where $\lambda$ is of the same form as stated in Theorem.  Since 
(6.3.2) is a $\fgl(V)$-module homomorphism, its kernel is 
$\mathop{\oplus}\limits_{\lambda: \lambda _{n+1}\geq 
n+2}W^*{}^{\lambda}$.  The fact that this kernel is generated by the 
elements of the least degreee is proved by the same arguments as for 
$\fosp(V)$.  \end{proof}

\section*{\S 7. The invariants of $\fspe(V)$}

First, let us describe certain tensor invariants. Let $\cT_1$ be
the set of matrices $A$ whose entries are equal to either 1 or 0,
with zeroes on the main diagonal and such that
$a_{ij}+a_{ji}=1$ for all offdiagonal entries. Set
$$
A_i=\sum a_{pq},\text{ where the sum runs over all the
elements strictly below the $i$-th row}.
$$
Define $|A|$ recursively: for $n=2$ set $|A|=0$ and for $n>2$ set
$$
|A|=|A^*|+\mathop{\sum}\limits_{i=1}^{n-2}a_{in}A^*_i+
\mathop{\sum}\limits_{1\leq
j<i<n}a_{in}a_{nj}+\mathop{\sum}\limits_{i>j}a_{ij}+\frac16n(n-1)(n-2),
$$
where $A^*$ is obtained from $A$ by striking out the last row
and the last column.

\ssbegin{7.1}{Lemma}
$Y^-=\mathop{\prod}\limits_{i<j}(E_{\bar i, j} - E_{\bar j,
i})=\mathop{\sum}\limits_{A\in\cT_1}(-1)^{|A|}E_A$, where the 
product runs over the lexicografically ordered set of pairs $i<j$
and $E_A=\prod E_{\bar i, j}^{a_{i, j}}$ and where the last
product is taken over the rows of the matrix
$A$ from left to right and downwards. \end{Lemma}

\begin{proof} Clearly, $Y^-$ is the product of $\frac 12n(n-1)$
factors. In each factor, select either $E_{\bar i, j}$ or $E_{\bar
j, i}$. In the first case, for $E_{\bar i, j}$, set $a_{i, j}=1$
and $a_{j, i}=0$ in the second case set the other way round.
We get a matrix with the properties desired. The sign is
obtained after reordering of the sequence of the $a_{ij}$:
$$
a_{12}a_{21}a_{13}a_{31}\dots a_{1n}a_{n1}\dots
a_{n-1,n}a_{n, n-1}\mapsto
a_{12}a_{13}\dots a_{1n}\dots a_{n-1,n}a_{n1}a_{n2}\dots
a_{n, n-1}.
$$
This is performed by induction: first, the pairs $a_{in}a_{ni}$
are moved to the end in increasing order, this accrues the
exponent of the sign with $\frac16n(n-1)(n-2)$,
then we reorder the elements with indices lesser
than $n$, which adds $|A^*|$, then the elements of
the sequence $a_{1n}a_{n1}a_{2n}a_{n2}\dots
a_{n-1,n}a_{n, n-1}$ are rearranged into the
sequence  $a_{1n}\dots
a_{n-1,n}a_{n1}a_{n2}\dots a_{n, n-1}$ which adds
$\mathop{\sum}\limits_{j<i}a_{in}a_{ni}$ to the exponent, and,
finally, the elements $a_{1n}$, \dots , $a_{n-1,n}$ are placed
onto the end of the $i$-th row adding
$\mathop{\sum}\limits_{i=1}^{n-2}a_{in}A^*_{i}$. Besides, if
$i>j$, then $E_{\bar i, j}$ enters $Y^-$ with a minus sign; this
adds $\mathop{\sum}\limits_{i>j}a_{ij}$. 
\end{proof}

\ssec{7.2} the numbers $i$ and $j$ we be referred to as {\it 
conjugate} if $i=\bar j$, i.e., they are equal but belong to copies of 
$\Nee$ of distinct ``parity''.  Let $\cT_2$ be the set of sequences of 
length $n^2$ considered as $n\times n$-tableaux filled in along 
columns and with the following properties: the numbers symmetric with 
respect to the main diagonal are conjugate, the $(i, j)$-th entry is 
occupied with one of the numbers $i$ or $\bar j$, the main diagonal is 
filled in with ``odd" numbers $\bar 1$, \dots , $\bar n$.  For every 
$L\in\cT_2$ determine the matrix $A=(a_{ij})$ by setting 
$a_{ij}=p(l_{ij})$, $n(L)=\#(\text{``even" elements in }L)$, 
$m(L)=m(A)=\mathop{\sum}\limits_{i\text{ is ``even"}}a_{ij}$, 
$\eps(L)=(-1)^{|A|+n(L)}$; let 
$m_k(L)=\frac{((n+k)!)^n}{(n+k-l_1)!\dots (n+k-l_n)!}$, where 
$l_i=\sum_j a_{ij}$.

\begin{Theorem} The elements 
$$
e_t\left (\mathop{\sum}\limits_{}(-1)^{km(L)}\eps(L)m_k(L)v^*_L\otimes 
v^*_{J_{k}}\right)\text{ for }L\in\cT_2 $$ and
$$
e_t\left (\mathop{\sum}\limits_{}\eps(L)m_0(L)v^*_L\otimes 
v^*_{I_{k}}\right )\text{ for }L\in\cT_2
$$
are
$\fspe(V)$-invariant. \end{Theorem}

\begin{proof} Let $r$ be an $(n+k)\times n$ rectangle filled in
along columns. Set $w=v^*_{J_{n+k}}$. Denote by
$w^{j_{1}}_{i_{1}}\dots w^{j_{p}}_{i_{p}}$ the tensor obtained
from $w$ by replacing the elements occupying positions
$i_{1}$, \dots , $i_{p}$ with numbers $j_{1}$, \dots , $j_{p}$,
respectively. Then
$$
e_r(E_{\bar n, j}w)=(-1)^{i-1}e_r(w^{j}_{i})(n+k), 
$$
where $i$ is any of 
the numbers of the positions occupied by $\bar n$.

If $E_{A_{n}}=\prod E_{\bar n, j}^{\alpha_{n, j}}$, the product
being ordered in order of increase of indices $j$, then 
$$
e_r(E_{A_{n}}w)=(-1)^{i_{1}-1+ \dots +
i_{l}-1}\frac{(n+k)!}{(n+k-l)!}e_r(w^{j_{1}, \dots , j_{l}}_{i_{1},
\dots , i_{l}}),
$$
where $\{j_{1}, \dots , j_{l}\}=\{j\mid \alpha_{n, j}\neq 0\}$
and $l=\mathop{\sum}\limits_{j}\alpha_{n, j}$, where $i_{1}<
\dots < i_{l}$.

Assume that $\{i_{1},
\dots , i_{l}\}=\{a+j_{1}-1,  \dots , a+j_{l}-1\}$, where $a$ is
the number of the first element in the $n$-th column of
tableau $r$. We thus get
$$
e_r(E_{A_{n}}w)=(-1)^{l\cdot a+ j_{1}+\dots +
j_{l}}\frac{(n+k)!}{(n+k-l)!}e_r(w^{j_{1}, \dots , j_{l}}_{i_{1},
\dots , i_{l}}).
$$
By continuing the process we get
$$
e_r(E_{A}w)=(-1)^{\eps(A)}\frac{[(n+k)!]^n}{(n+k-l_{1})!\dots
(n+k-l_{n})!}e_r(v_{I_{A}}^*), 
$$
where $\eps(A)=a_1+\dots +a_n+n(A_{n})+n(A_{n-1})+\dots +n(A_{1})$ and 
$a_{i}$ is the number of the first element in the $i$th column and 
where $n(A_{i})$ is equal to the sum of the numbers of the places 
occupied by the 1's, and where $I_{A}$ coinsides with $J_{n+k}$ everywhere 
unless the $a_{ij}=1$; then the $(ij)$-th entry of $I_{A}$ is occupied 
by $j$.

If $I$ and $J$ are two sequences and $t$ and $s$ are two
tableaux of the same form such that after filling $t$ with the
elements from $I$ and $s$ with the
elements from $J$ one gets geometrically identical pictures,
then $\sigma t=s$ implies $\sigma I=J$. Indeed,
$$
I(\sigma^{-1}\alpha)=t(\sigma^{-1}\alpha)=s(\alpha)=J(\alpha).
$$
Therefore, if $\sigma t=r$, we have
$$
\renewcommand{\arraystretch}{1.4}
\begin{array}{l}
    Y^-e_t(v^*_{J_{n}}v^*_{J_{k}})=
Y^-e_{\sigma^{-1}r}(v^*_{\sigma^{-1}(J_{n+k})})=\\
Y^-\sigma^{-1}e_{r\sigma}(\sigma^{-1}v^*_{J_{n+k}})\cdot
c(J_{n+k}, \sigma)=
c(J_{n+k}, \sigma)\sigma^{-1}Y^-e_{r}(v^*_{J_{n+k}})
\end{array}
$$
because thanks to the fact that $J_{n+k}$ only contains ``odd"
elements $c(J_{n+k}, \sigma)=\sign (\sigma)$.

Hence,
$$
\renewcommand{\arraystretch}{1.4}
\begin{array}{l}
Y^-e_t(v^*_{J_{n}}v^*_{J_{k}})= \sign(\sigma)\sigma^{-1}e_r\left 
(\mathop{\sum}\limits_A(-1)^{\eps(A)} 
\frac{[(n+k)!]^n}{(n+k-l_{1})!\dots (n+k-l_{n})!} v_{I_{A}}^*\right )=\\
\sign(\sigma)e_t\left (\mathop{\sum}\limits_A(-1)^{\eps(A)} 
\frac{[(n+k)!]^n}{(n+k-l_{1})!\dots (n+k-l_{n})!} 
\sigma^{-1}v^*_{I_{A}}\right )= \\
\sign(\sigma)e_t\left 
(\mathop{\sum}\limits_A(-1)^{\eps(A)} 
\frac{[(n+k)!]^n}{(n+k-l_{1})!\dots (n+k-l_{n})!} c(I_{A}, 
\sigma)v^*_{J_{A}}\otimes v_{J_{k}}\right ).
\end{array}
$$
where
$$
c(I_{A}, \sigma)=|A_{2}|\cdot k+|A_{4}|\cdot k\dots =
k\mathop{\sum}\limits_{i \text{ is even}}a_{ij},
$$
and where $J_{A}$ coinsides with $J_{n}$ everywhere unless where 
$a_{ij}=1$, then the $(i,j)$the position is occupied by $j$.  We are 
done.
\end{proof}

\ssbegin{7.3}{Theorem} The algebra of $\fspe(V)$-invariant
polynomials is generated by the following elements

{\em i)} $(v_\alpha, v_\beta)$ for $\alpha, \beta\in R_W$;

{\em ii)} $PPf_{k}(J)=\mathop{\sum}\limits_{} 
(-1)^{(k-1)m(L)}\eps(L)m_{k-1}(L)P_t(L*J_k, J)$ for any $t$-standard 
sequence $J$;

{\em iii)}
$PPf_{-k}(J)=\mathop{\sum}\limits_{}\eps(L)m_{0}(L)P_t(L*I_{k+1},
J)$ for any $s$-standard sequence $J$,

\noindent where $k\geq 1$ and sums run over $l\in\cT_2$.
\end{Theorem}

\begin{proof} is similar to that of 5.2.
\end{proof}

\end{document}